\documentclass[11pt,reqno]{amsart}
\usepackage{amsfonts}
\usepackage{amsmath,amsthm,amsfonts,amssymb}
\newtheorem{lemma}{Lemma}
\newtheorem{theorem}{Theorem}
\newtheorem{corollary}{Corollary}

\def\bl{\begin{lemma}}
\def\bt{\begin{theorem}}
\def\el{\end{lemma}}
\def\et{\end{theorem}}
\def\bp{\begin{proof}}
\def\ep{\end{proof}}
\def\bc{\begin{corollary}}
\def\ec{\end{corollary}}

\def\mb{\mathbb}

\def\O{\Omega}
\def\a{\alpha}
\def\b{\beta}

\def\s{\sigma}

\def\-{\setminus}
\def\s{\sigma}

\def\ov{\overline}
\def\lt{\left}
\def\rt{\right}
\def\+{\bigcup}
\def\.{\bigcap}

\def\ll{\langle}
\def\rl{\rangle}

\def\te{\theta}

\title[A Schwarz-Pick lemma]{A Schwarz-Pick lemma for the modulus of holomorphic mappings from the polydisk into the
unit ball}
\thanks{Research supported by the National Natural Science Foundation of China
(No. 11201199) and by the Scientific Research Foundation of Jinling
Institute of Technology (No. Jit-b-201221).}

\author {Shaoyu Dai and Yifei Pan}

\address{Department of Mathematics, Jinling Institute of
Technology, Nanjing 211169, China}
\address{\it E-mail address: dymdsy@163.com}

\address{School of Mathematics and Informatics, Jiangxi Normal University,
Nanchang 330022, China}
\address{Department of Mathematical Sciences, Indiana University -
Purdue University Fort Wayne, Fort Wayne, IN 46805-1499, USA}
\address{\it E-mail address: pan@ipfw.edu}

\begin{document}

\numberwithin{equation}{section}

\begin{abstract}
In this paper we prove a Schwarz-Pick lemma for the modulus of
holomorphic mappings from the polydisk into the unit ball. This
result extends some related results.
\end{abstract}

\maketitle



\section{Introduction}

Let $\mathbb{D}$ be the unit disk in $\mathbb{C}$, $\mathbb{D}^n$
and $\mathbb{B}_n$ be the polydisk and the unit ball in
$\mathbb{C}^n$ respectively. For $z=(z_1,\cdots,z_n)$ and
$z'=(z'_1,\cdots,z'_n)\in\mb C^n$, denote $\lt\ll z,z'\rt\rl=z_1\ov
z_1'+\cdots+z_n\ov z_n'$ and $|z|=\lt\ll z,z\rt\rl^{1/2}$. Let
$\O_{X,Y}$ be the class of all holomorphic mappings $f$ from $X$
into $Y$, where $X$ is a domain in $\mathbb{C}^n$ and $Y$ is a
domain in $\mathbb{C}^m$. For $f\in\O_{X,Y}$ and $j=1,\cdots,n$,
define
\begin{equation}\label{1}
|\nabla|f|(z)|=\sup_{\b\in
\mathbb{C}^n,\,|\b|=1}\lt(\lim_{t\in\mathbb{R},\,t\rightarrow0^+}\frac{|f|(z+t\b)-|f|(z)}{t}\rt),\
\ \ \ \ \ z\in X;
\end{equation}
\begin{equation}\label{2}
|\nabla_j|f|(z)|=\sup_{\b\in
\mathbb{C},\,|\b|=1}\lt(\lim_{t\in\mathbb{R},\,t\rightarrow0^+}\frac{|f|(z_1,\cdots,z_{j-1},z_j+t\b,z_{j+1},\cdots,z_n)-|f|(z)}{t}\rt),\
\ \ \ \ \ z\in X,
\end{equation}
where $f=(f_1,\cdots,f_m)$,
$|f|=(|f_1|^2+\cdots+|f_m|^2)^{\frac{1}{2}}$ and
$z=(z_1,\cdots,z_n)$. Some calculation for $|\nabla|f||$ and
$|\nabla_j|f||$ will be given in Section 2.

For $f\in\O_{\mathbb{D},\mathbb{D}}$, the classical Schwarz-Pick
lemma says that
\begin{equation}\label{3}
|f'(z)|\leq\frac{1-|f(z)|^2}{1-|z|^2},\ \ \ \ \ \ z\in\mathbb{D}.
\end{equation}
This inequality does not hold for $f\in\O_{\mathbb{D},\mathbb{B}_m}$
with $m\geq2$. For instance, the mapping
$f(z)=\frac{1}{\sqrt{2}}(z,1)$ satisfies
$$|f'(0)|=\sqrt{1-|f(0)|^2}>1-|f(0)|^2.$$
However Pavlovi$\acute{c}$ \cite{P} found that \eqref{3} can also be
written as
\begin{equation}\label{4}
|\nabla|f|(z)|\leq\frac{1-|f(z)|^2}{1-|z|^2},\ \ \ \ \ \
z\in\mathbb{D},
\end{equation}
since \eqref{11}. In \cite{P}, Pavlovi$\acute{c}$ proved that this
form \eqref{4} can be extended to $\O_{\mathbb{D},\mathbb{B}_m}$ and
obtained the same inequality for $f\in\O_{\mathbb{D},\mathbb{B}_m}$.
Recently, we \cite{D3} proved that the form \eqref{4} also can be
extended to $\O_{\mathbb{B}_n,\mathbb{B}_m}$ and obtained the
following inequality for $f\in\O_{\mathbb{B}_n,\mathbb{B}_m}$:
\begin{equation}\label{5}
|\nabla|f|(z)|\leq\frac{1-|f(z)|^2}{1-|z|^2},\ \ \ \ \ \
z\in\mathbb{B}_n.
\end{equation}
In view of the above results, it is interesting for us to consider
that if there are some similar results for
$f\in\O_{\mathbb{D}^n,\mathbb{B}_m}$.

For $f\in\O_{\mathbb{D}^n,\mathbb{D}}$, it is well known \cite{R,K}
that
\begin{equation}\label{6}
\sum^n_{j=1}(1-|z_j|^2)|f'_{z_j}(z)|\leq 1-|f(z)|^2
\end{equation}
for any $z=(z_1,\cdots,z_n)\in\mathbb{D}^n$. This inequality does
not hold for $f\in\O_{\mathbb{D}^n,\mathbb{B}_m}$ with $m\geq2$. For
instance, the mapping
$f(z)=\frac{1}{\sqrt{3}}(z_1,z_2+0.1)\in\O_{\mathbb{D}^2,\mathbb{B}_2}$
satisfies
$$\sum^2_{j=1}|f'_{z_j}(0)|=\frac{2}{\sqrt{3}}>1-|f(0)|^2.$$
Similarly to \eqref{4}, we find that \eqref{6} can be written as
\begin{equation}\label{7}
\sum^n_{j=1}(1-|z_j|^2)|\nabla_j|f|(z)|\leq 1-|f(z)|^2
\end{equation}
for any $z=(z_1,\cdots,z_n)\in\mathbb{D}^n$, since \eqref{18}. In
view of \eqref{4} and \eqref{5}, the obvious question is that if the
form \eqref{7} can be extended to $\O_{\mathbb{D}^n,\mathbb{B}_m}$
with $m\geq2$. The following example shows that the form \eqref{7}
can not completely be extended to $\O_{\mathbb{D}^n,\mathbb{B}_m}$
with $m\geq2$: the mapping
$f(z)=\frac{1}{\sqrt{2}}(z_1,z_2)\in\O_{\mathbb{D}^2,\mathbb{B}_2}$
satisfies
$$\sum^2_{j=1}|\nabla_j|f|(0)|=\sqrt{2}>1-|f(0)|^2,$$
since $f(0)=0$ and $|\nabla_j|f|(0)|=|f'_{z_j}(0)|$ for $j=1,2$ by
\eqref{12}. However we find that the form \eqref{7} holds for
$f\in\O_{\mathbb{D}^n,\mathbb{B}_m}$ at the point $z\in\mathbb{D}^n$
with $f(z)\neq0$. Precisely:

\bt\label{th1} Let $f: \mathbb{D}^n\rightarrow\mathbb{B}_m$ be a
holomorphic mapping with $m\geq2$. Then
\begin{equation}\label{8}
\sum^n_{j=1}(1-|z_j|^2)|\nabla_j|f|(z)|\leq 1-|f(z)|^2, \ \ \ \
\mbox{if}\quad f(z)\neq0
\end{equation}
and
\begin{equation}\label{28}
\sum^n_{j=1}(1-|z_j|^2)^2|\nabla_j|f|(z)|^2\leq 1, \ \ \ \
\mbox{if}\quad f(z)=0
\end{equation}
for any $z=(z_1,\cdots,z_n)\in\mathbb{D}^n$. \et

The above theorem is the main result in this paper. Note that the
inequality in \eqref{28} always holds whether if $f(z)=0$ or
$f(z)\neq0$. When $f(z)\neq0$, there is a better inequality, which
is \eqref{8}. Theorem 1 is coincident with \eqref{5} when $n=1$. In
addition, \eqref{8} and \eqref{28} are sharp. For example, the
mapping $f(z)=\frac{1}{\sqrt{2}}
\lt(\frac{\frac{1}{2}-z_1}{1-\frac{1}{2}z_1},
\frac{\frac{1}{2}-z_2}{1-\frac{1}{2}z_2}\rt)\in\O_{\mathbb{D}^2,\mathbb{B}_2}$
satisfies the equality in \eqref{8} at $z=0$; the mapping
$f(z)=\frac{1}{\sqrt{2}}(z_1,z_2)\in\O_{\mathbb{D}^2,\mathbb{B}_2}$
satisfies the equality in \eqref{28} at $z=0$.

In Section 2, some calculation for $|\nabla|f||$ and $|\nabla_j|f||$
will be given. In Section 3, we will give the proof of Theorem
\ref{th1} and some remarks for the equality cases in Theorem 1.

\section{Some calculation for
$|\nabla|f||$ and $|\nabla_j|f||$}

For $f\in\O_{X,Y}$ with $X\subset\mathbb{C}^n$ and
$Y\subset\mathbb{C}^m$, by \eqref{1} we know that if $|f|(z)\neq0$
then $|f|$ is $\mathbb{R}$-differentiable at $z$ and $\nabla |f|$ is
the ordinary gradient;  if $|f|(z)=0$ then $|f|$ is not
$\mathbb{R}$-differentiable at $z$ and $\nabla |f|$ is not the
ordinary gradient. From Section 2 in \cite{D3}, we have the
following \eqref{9}-\eqref{11}. For $f\in\O_{X,Y}$,
\begin{equation}\label{9}
|\nabla|f|(z)|=
\begin{cases}
\frac{1}{|f(z)|}\lt|\lt(\lt\ll f'_{z_1}(z),
f(z)\rt\rl,\cdots,\lt\ll f'_{z_n}(z), f(z)\rt\rl\rt)\rt|, &\mbox{if}\quad f(z)\neq0;\\
\underset{\b\in \mathbb{C}^n,\,|\b|=1}{\sup}|Df(z)\cdot \b|,
&\mbox{if}\quad f(z)=0,
\end{cases}
\end{equation}
where $z=(z_1,\cdots,z_n)\in X$ and $Df(z)\cdot \b$ is the Fr\'echet
derivative of $f$ at $z$ in the direction $\b$. Then for
$f\in\O_{X,Y}$ with $X\subset\mathbb{C}$,
\begin{equation}\label{10}
|\nabla|f|(z)|=
\begin{cases}
\frac{1}{|f(z)|}\lt|\lt\ll f'(z),
f(z)\rt\rl\rt|, &\mbox{if}\quad f(z)\neq0;\\
|f'(z)|, &\mbox{if}\quad f(z)=0.
\end{cases}
\end{equation}
In particular, for $f\in\O_{X,Y}$ with $X\subset\mathbb{C}$ and
$Y\subset\mathbb{C}$,
\begin{equation}\label{11}
|\nabla|f|(z)|=|f'(z)|.
\end{equation}

Then by \eqref{2} and \eqref{10}, we get that for $f\in\O_{X,Y}$ and
$j=1,\cdots,n$,
\begin{equation}\label{12}
|\nabla_j|f|(z)|=
\begin{cases}
\frac{1}{|f(z)|}\lt|\lt\ll f'_{z_j}(z),
f(z)\rt\rl\rt|, &\mbox{if}\quad f(z)\neq0;\\
|f'_{z_j}(z)|, &\mbox{if}\quad f(z)=0,
\end{cases}
\end{equation}
where $z=(z_1,\cdots,z_n)\in X$. Note that for the case that
$f(z)\neq0$, if $f'_{z_j}(z)$ and $f(z)$ are collinear, then
$|\nabla_j|f|(z)|=|f'_{z_j}(z)|$; if not, then
$|\nabla_j|f|(z)|\neq|f'_{z_j}(z)|$. In particular, for
$f\in\O_{X,Y}$ with $Y\subset\mathbb{C}$,
\begin{equation}\label{18}
|\nabla_j|f|(z)|=|f'_{z_j}(z)|.
\end{equation}

\section{Proof of Theorem \ref{th1}}

First we give one lemma.

\bl\label{l1}  Let
$f(z)=\underset{\a}{\sum}a_{\a}z^{\a}\in\O_{\mathbb{D}^n,\mathbb{B}_m}$,
where $z=(z_{1},\cdots,z_{n})$, $\a=(\a_{1},\cdots,\a_{n})$,
$z^\a=z_1^{\a_1},\cdots,z_n^{\a_n}$, $f=(f_1,\cdots,f_m)$,
$f_j(z)=\sum\limits_{\a}a_{j,\a}z^{\a}$ and
$a_\a=(a_{1,\a},\cdots,a_{m,\a})$. Then
\begin{equation}\label{30}
\sum_\a|a_\a|^2\leq1.
\end{equation}

 \el

\bp  For $0<\s<1$, we have
\begin{equation*}
\begin{split}
1&\ge\frac1{(2\pi)^n}\int_0^{2\pi}\!\!\!\!\cdots\!\int_0^{2\pi}
|f(\s e^{i\te_1},\cdots,\s e^{i\te_n})|^2d\te_1\cdots d\te_n\\
&=\frac1{(2\pi)^n}\sum_{j=1}^m\int_0^{2\pi}\!\!\!\!\cdots\!\int_0^{2\pi}
|f_j(\s e^{i\te_1},\cdots,\s e^{i\te_n})|^2d\te_1\cdots
d\te_n\\
&=\sum_{j=1}^m\sum_\a|a_{j,\a}|^2\s^{2|\a|}\\
&=\sum_\a|a_\a|^2\s^{2|\a|},
\end{split}
\end{equation*}
where $|\a|=\overset{n}{\underset{j=1}{\sum}}\a_{j}$. Letting
$\s\to1$ gives \eqref{30}.  \ep

Now we give the proof of Theorem \ref{th1}.

\noindent{\it Proof of Theorem \ref{th1}.}\quad First we prove the
case that $z=0$.

Therefore we need to prove that
\begin{equation}\label{13}
\begin{cases}
\sum^n_{j=1}|\nabla_j|f|(0)|\leq 1-|f(0)|^2, &\mbox{if}\quad f(0)\neq0;\\
\sum^n_{j=1}|\nabla_j|f|(0)|^2\leq 1, &\mbox{if}\quad f(0)=0.
\end{cases}
\end{equation}
By \eqref{12}, it suffices to prove that
\begin{equation}\label{14}
\sum^n_{j=1}\lt|\lt\ll f'_{z_j}(0),\frac{f(0)}{|f(0)|}
\rt\rl\rt|\leq 1-|f(0)|^2, \ \ \mbox{if}\quad f(0)\neq0
\end{equation}
and
\begin{equation}\label{15}
\sum^n_{j=1}|f'_{z_j}(0)|^2\leq 1, \ \ \mbox{if}\quad f(0)=0.
\end{equation}
Obviously, \eqref{15} holds by Lemma \ref{l1}. For \eqref{14}, let
$$h(z)=\lt\ll f(z),\frac{f(0)}{|f(0)|}\rt\rl,\ \ \ \ \ \
z\in\mathbb{D}^n.$$ Then $h(z)$ is a holomorphic function from
$\mathbb{D}^n$ into $\mathbb{D}$, $h(0)=|f(0)|$, and for
$j=1,\cdots,n$,
\begin{equation}\label{16}
h'_{z_j}(0)=\lt\ll f'_{z_j}(0),\frac{f(0)}{|f(0)|}\rt\rl,
\end{equation}
where $z=(z_1,\cdots,z_n)$. Applying \eqref{6} to $h$ and by
\eqref{16} we get
\begin{equation*}
\begin{split}
\sum^n_{j=1}\lt|\lt\ll f'_{z_j}(0),\frac{f(0)}{|f(0)|}
\rt\rl\rt|&=\sum^n_{j=1}|h'_{z_j}(0)|\\
&\leq 1-|h(0)|^2\\
&=1-|f(0)|^2.
\end{split}
\end{equation*}
Then \eqref{14} is proved. Therefore \eqref{13} is proved.

Now we prove the case that $z=p\neq0$.

Let $p=(p_1,\cdots,p_n)$ and
$$g(w)=f(\varphi(w)),\ \ \ \ \ \
w=(w_1,\cdots,w_n)\in\mathbb{D}^n,$$ where
$\varphi(w)=(\varphi_1(w_1),\cdots,\varphi_n(w_n))$,
$\varphi_j(w_j)=\frac{p_j-w_j}{1-\overline{p_j}w_j}$ for
$j=1,\cdots,n$. Then $g(w)$ is a holomorphic mapping from
$\mathbb{D}^n$ into $\mathbb{B}_m$, $g(0)=f(p)$, and for
$j=1,\cdots,n$,
\begin{equation}\label{17}
g'_{w_j}(0)=f'_{z_j}(p)(-1+|p_j|^2).
\end{equation}
For the case that $f(p)\neq0$, applying \eqref{14} to $g$ and by
\eqref{12}, \eqref{17} we get
\begin{equation*}
\begin{split}
\sum^n_{j=1}(1-|p_j|^2)|\nabla_j|f|(p)|&=\sum^n_{j=1}(1-|p_j|^2)\lt|\lt\ll
f'_{z_j}(p),\frac{f(p)}{|f(p)|} \rt\rl\rt|\\&=\sum^n_{j=1}\lt|\lt\ll
g'_{w_j}(0),\frac{g(0)}{|g(0)|} \rt\rl\rt|\\&\leq
1-|g(0)|^2\\&=1-|f(p)|^2.
\end{split}
\end{equation*}
For the case that $f(p)=0$, applying \eqref{15} to $g$ and by
\eqref{12}, \eqref{17} we get
\begin{equation*}
\begin{split}
\sum^n_{j=1}(1-|p_j|^2)^2|\nabla_j|f|(p)|^2&=\sum^n_{j=1}(1-|p_j|^2)^2|f'_{z_j}(p)|^2\\&=\sum^n_{j=1}|g'_{w_j}(0)|^2\\&\leq
1.
\end{split}
\end{equation*}
Then the theorem is proved.
 \qed

In the following, we give some remarks for the equality cases in
Theorem 1.

{\bf Remark 1.} When $n=1$, \eqref{8} and \eqref{28} reduce to
\eqref{5}. The equality case in \eqref{5} has been discussed in
\cite{D3}.

{\bf Remark 2.} When $n\geq2$, if the equality in \eqref{28} holds
at some point $p=(p_1,\cdots,p_n)$, then the structure of the
expression of $f$ will be controlled. Precisely:
$$f(z)=\sum^n_{j=1}f'_{z_j}(p)(-1+|p_j|^2)\frac{p_j-z_j}{1-\overline{p_j}z_j},\ \ \ \ \ \
z\in\mathbb{D}^n,$$ which is obvious by the proof of Theorem 1,
Lemma 1 and \eqref{12}.

{\bf Remark 3.} When $n\geq2$, if the equality in \eqref{8} holds at
some point $p=(p_1,\cdots,p_n)$, then the following discussion shows
that the equality at $p$ is not enough to control the structure of
the expression of $f$. By the proof of Theorem 1, we know that the
key to the extremal problem of \eqref{8} at the point $p$ is to
solve the extremal problem of \eqref{6} at $z=0$. That is: for
$h\in\O_{\mathbb{D}^n,\mathbb{D}}$, if
$\sum^n_{j=1}|h'_{z_j}(0)|=1-|h(0)|^2$, then what the structure of
the expression of $h$ is. By the proof of \eqref{6} in \cite{R}, we
only need to consider this problem: for
$h\in\O_{\mathbb{D}^n,\mathbb{D}}$ with $h(0)=0$, if
$\sum^n_{j=1}|h'_{z_j}(0)|=1$, then what the structure of the
expression of $h$ is. However, the following examples show that the
condition $\sum^n_{j=1}|h'_{z_j}(0)|=1$ can not control the higher
order terms in the expansion of $h$. Consequently, the structure of
the expression of $h$ can not be controlled.\\
Examples:
$$g(z)=\frac{1}{2}z_1+\frac{1}{2}z_2\in\O_{\mathbb{D}^2,\mathbb{D}};$$
$$\tilde{g}(z)=\frac{\frac{1}{2}z_1+\frac{1}{2}z_2-z_1z_2}{1-\frac{1}{2}z_1-\frac{1}{2}z_2}
\in\O_{\mathbb{D}^2,\mathbb{D}}.$$ Although the above two functions
satisfy $g(0)=\tilde{g}(0)=0$, $g'_{z_1}(0)=\tilde{g}'_{z_1}(0)$,
$g'_{z_2}(0)=\tilde{g}'_{z_2}(0)$ and
$\sum^n_{j=1}|g'_{z_j}(0)|=\sum^n_{j=1}|\tilde{g}'_{z_j}(0)|=1$, the
expression of $g$ has no higher order terms and the expression of
$\tilde{g}$ has some higher order terms.


\end{document}